\documentclass[10pt]{article}
\advance \topmargin by -\headheight
\advance \topmargin by -\headsep
\evensidemargin \oddsidemargin
\usepackage{epsfig}
\usepackage{amsthm}
\usepackage{amssymb}
\usepackage{amsmath}

\theoremstyle{plain}
\newtheorem{Lemma}{Lemma}
\newtheorem*{Theorem}{Theorem}
\theoremstyle{definition}

\newcommand{\A}{{\mathcal A}}
\newcommand{\B}{{\mathcal B}}

\let\:=\colon

\begin{document}

\title{A Large Dihedral Symmetry of the Set of Alternating Sign Matrices}
\author{{\sc Benjamin Wieland}\thanks{
The author is grateful to James Propp for discussions and comments on 
the paper.
}\\ \texttt{wieland@math.uchicago.edu}}
\date{June 30, 2000}

\maketitle

\begin{abstract}
We prove a conjecture of Cohn and Propp, which refines a conjecture of
Bosley and Fidkowski about the symmetry of the set of alternating sign 
matrices (ASMs).  We examine data arising from the representation of 
an ASM as a collection of paths connecting $2n$ vertices and show it 
to be invariant under the dihedral group $D_{2n}$ rearranging those 
vertices, which is much bigger than the group of symmetries of the 
square.  We also generalize conjectures of Propp and Wilson relating 
some of this data for different values of $n$.
\end{abstract}

\section{Introduction}
\label{intro}
\begin{figure}[htb]
\begin{center}
\epsfig{file=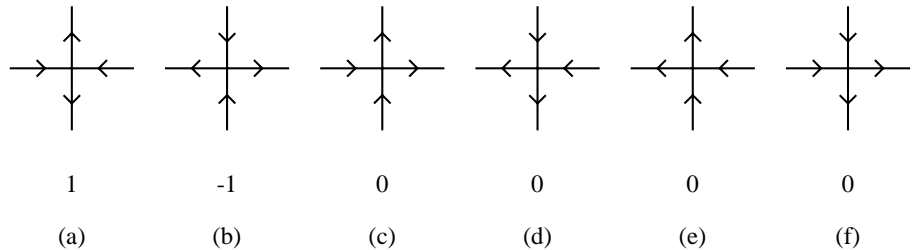,height=1.3in}
\end{center}
\caption{The six vertex configurations for the square ice model.}
\label{ice}
\end{figure}

In statistical mechanics, the {\em square ice}\/ model represents an 
ice crystal as a directed graph with 
oxygen atoms at the vertices and hydrogen atoms on the edges.  A 
hydrogen atom is shared between two oxygen atoms, but it is covalently 
bonded to one atom, and it is hydrogen bonded to the other; correspondingly, 
an edge is incident to two vertices, but favors one by pointing at it.  
Each oxygen has two covalently bonded hydrogens, so each vertex has 
two edges pointing in and two pointing out, as in Figure~\ref{ice}.  
The model is called 
``square'' because the graph is taken to be part of the square grid.

\begin{figure}[htb]
\begin{center}
\epsfig{file=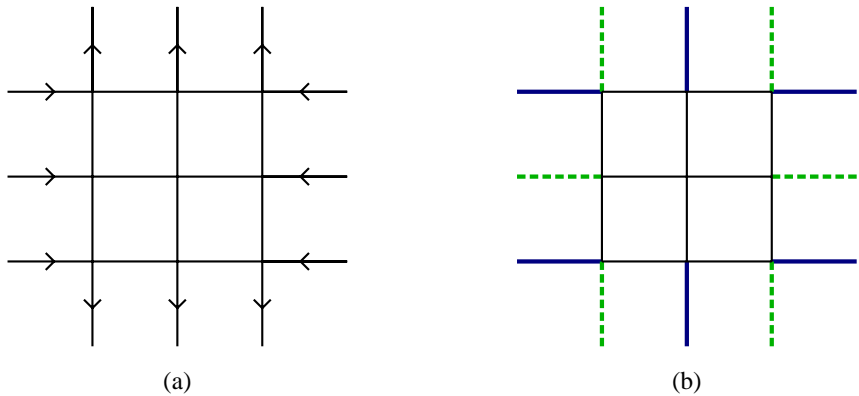,height=2in}
\end{center}
\caption{The graph and boundary conditions to put square ice
configurations in bijection with ASMs for $n=3$ and the corresponding
boundary conditions for paths.}
\label{boundary}
\end{figure}

In statistical mechanics, the model is interesting on a torus, or on a
finite square with unrestricted boundary conditions, but
combinatorialists are most interested in applying these restrictions
to a graph with a particular set of boundary conditions, which will
not allow it to wrap around a torus.  At the boundary where the square
ends, there are vertices with only one edge.  We require that these
edges point into the square if they are horizontal, and out of the
square if they are vertical, as in Figure~\ref{boundary}(a).  
With these boundary conditions, orientations of a finite square grid
are in bijection
with a certain set of  matrices, called {\em alternating sign
matrices}\/ (ASMs), which are defined in Section~\ref{asms and
heights}.  We identify the graphs with the matrices, and
refer to them also as ASMs to indicate that they have the proper
boundary conditions.  ASMs are of greater interest in combinatorics;
for example the number of ASMs was proved in \cite{Z,K} to be
$\prod_{i=0}^{n-1}\frac{(3i+1)!}{(n+i)!}$, settling a long-standing
conjecture of \cite{MRR82}. See \cite{B} for a history of the problem.

If we color the vertices alternately black and white, and, if we color 
the edges blue and green depending on the color of the initial 
vertices, then the rule of two edges in and two out becomes the 
requirement that each vertex lie on two blue edges and two green 
edges.  As the old boundary conditions determined the orientation of 
each edge on the boundary, the new one determine its color.  
Specifically, the boundary edges alternate in color, as in 
Figure~\ref{boundary}(b).

Because no vertex lies on more than two blue edges, we may start at a 
blue edge on the boundary and follow that edge into the graph.  At 
each interior vertex, there are two blue edges; so there is only one 
way to leave the vertex while remaining on blue and not retracing the 
path.  Eventually this path ends somewhere else along the boundary. 
Thus, the path associates the initial vertex with the final vertex.
This applies to each boundary vertex lying on a blue edge, so these 
vertices are divided into pairs by the blue paths of the ASM.

By focusing on the boundary, we may ignore the square structure of the 
graph, and think only of the outer ring of vertices in their cyclic 
order.  Although a square must be rotated by multiples of $\frac\pi2$, 
we may rotate the boundary a much smaller amount, and ask how many 
ASMs produce the new pairing.  The central result of this paper is 
that there are as many ASMs with the new pairing as with the old one, 
as conjectured by Carl Bosley and Lukasz Fidkowski\cite{BF}; 
in fact, we construct an explicit bijection: from each 
ASM, we construct a new one, with the blue pairing rotated one step 
clockwise.  Absolute position does not matter: a pairing may have 
connected two corners of the square, but now connects vertices in the 
middles of sides.

By considering the green edges and the closed loops, Henry Cohn and 
James Propp \cite{CP} refined this conjecture to the 
form that we prove here.  The boundary edges alternate blue and green.  
Exactly as with the blue ones, the green edges start paths joining the 
other half of the boundary vertices.  The construction causes the 
green pairing to rotate not clockwise, but counterclockwise.  Because 
of these opposing rotations, we call the bijection {\em gyration}.  If 
we start at an interior vertex, then there are two blue edges, which 
may be followed to produce two blue paths.  Either they both terminate 
on the boundary, so have been considered in the discussion of 
pairings, or they rejoin producing a closed loop.  While gyration 
turns some loops from one color to the other, the total number remains 
constant.

Section~\ref{asms and heights} defines ASMs and another equivalent
representation, called a height function, for which gyration is
particularly simple.  Section~\ref{theorem} contains the statement of
our theorem on the properties of gyration.  In Section~\ref{gyration},
we define gyration.  In Section~\ref{gyration2}, we split gyration
into two steps, which are useful because they have similar properties
to gyration itself, but which violate the boundary conditions.
Section~\ref{lemma} proves a lemma, which identifies monochromatic
paths before and after the these half-steps of gyration.  It is the
key step.  Section~\ref{proof} combines the lemma with the boundary
conditions to prove the theorem.  In Section~\ref{generalizations} we
describe other contexts for gyration and state some open problems.

\begin{figure}[htb]
\begin{center}
\epsfig{file=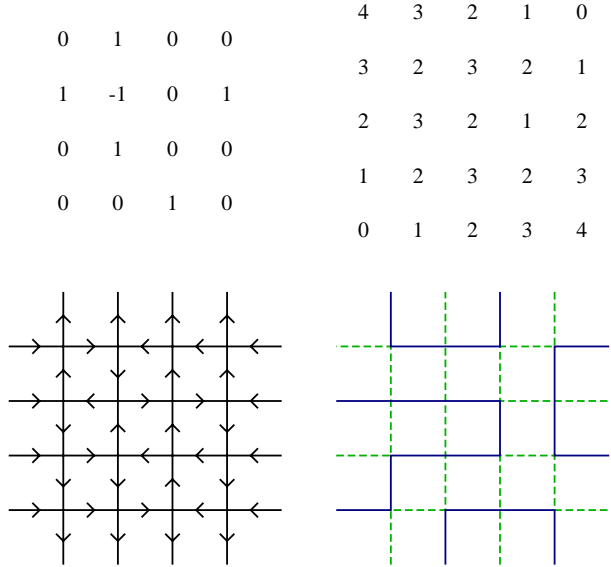,height=3.0in}
\end{center}
\caption{An ASM represented as a matrix, 
a square ice configuration, and a colored graph.}
\label{same}
\end{figure}

\section{Relation to ASMs and Height Functions}
\label{asms and heights}
An {\em alternating sign matrix} (ASM) of order $n$ is an $n\times n$
square matrix such that each entry is one of $1$, $-1$, and $0$, and
the nonzero entries in each row and column alternate in sign and sum
to $1$. Let $\A_n$ denote the set ASMs of order $n$.  ASMs are in
bijection with another class of matrix, called the corner sum matrix
in \cite{RR} and the skewed summation in \cite{EKLP}.  These are
$(n+1)\times(n+1)$ matrices with integral entries such that
horizontally and vertically adjacent entries differing by $1$ and the
entries on the edge fixed at the values shown in
Figure~\ref{same}. See \cite{RR} or \cite{EKLP} for the specific
bijection.  We adopt the statistical mechanics term {\em height
function}\/ for this type of matrix. 
We will not make use of height functions, but we mention them
because the 
simplest definition of gyration is in terms of height functions (see
Section~\ref{gyration}), and this definition can apply to many
height functions. The entries of a height function matrix for an ASM
should be thought of as ofset from the entries of the ASM
itself. If the entries of the ASM lie on vertices in the square grid,
then the entries of the height function lie on the faces of that graph.

We may choose the subgraph and some boundary conditions so that the
set of configurations of square ice 
is in bijection with ASMs of order $n$ (see,
e.g., \cite{EKLP}).  Take the $n^2$ vertices with both coordinates
between $1$ and $n$, which are called {\em interior vertices}.  Take
all edges incident to these vertices; doing so requires that we take
the $4n$ additional neighboring vertices, which have one coordinate
from $1$ to $n$ and the other either $0$ or $n+1$; call these vertices
{\em endpoints}. Denote this graph $L_n$.  The interior vertices have
degree $4$, and must be in one of the six vertex configurations,
whereas the endpoints have degree $1$ and cannot.  Instead, we require
that the vertical edges incident to endpoints be directed out of the
square and that the horizontal edges be directed in, as in
Figure~\ref{boundary}(a).

To turn a square ice configuration with such boundary conditions into 
an ASM, replace the vertical-out, horizontal-in vertices with $1$s, 
the vertical-in, horizontal-out vertices with $-1$s, and all other 
vertices with $0$s, as indicated in Figure~\ref{ice}.  Since the four 
configurations that become $0$s have two horizontal edges oriented in 
the same absolute direction, the only way a line of horizontal edges, 
traversed from left to right, changes from edges pointing right to 
edges pointing left is by passing through a $1$; similarly, a change back 
comes from passing through a $-1$.  Any vertex at which they do not 
change must come from a $0$.  Consider two vertices in the same row 
that become nonzero and have only $0$s between them.  All the edges 
between them must point in the same direction.  The nonzero vertex in 
that direction is a $1$, and the other a $-1$.  Thus the square ice 
model produces matrices with alternating nonzero entries.  The 
boundary conditions require that the outermost nonzero entries both be 
$1$s; so the sum of the entries in a row is $1$.  The vertical 
situation is similar, though the role of ``in'' and ``out'' switches, 
both in the configuration that becomes $1$ and in the boundary conditions.

\begin{figure}[htb]
\begin{center}
\epsfig{file=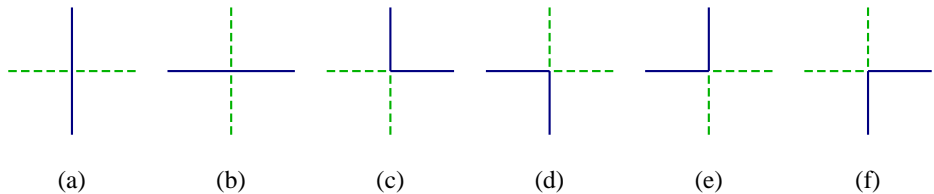,height=1in}
\end{center}
\caption{The six vertex possibilities for vertices in the colored 
graph version of square ice.}
\label{2factor}
\end{figure}

If we give a vertex $(x,y)$ the same parity as $x+y$, and if we color 
blue those edges directed from odd to even and green those directed 
from even to odd, then we get a graph in which every interior vertex 
has two incident green edges and two incident blue edges, except the 
endpoints, which alternate in the color of their incident edges.  In 
figures, we use solid and dashed lines to represent blue and 
green, respectively.  Let us give an endpoint the same color as its 
incident edge.  The six types of vertices in Figure~\ref{ice} turn 
into the six types of vertices in Figure~\ref{2factor} in the same 
order if the vertex is odd.  If the vertex is even, then \ref{ice}(a) 
and \ref{ice}(b) become \ref{2factor}(b) and \ref{2factor}(a), 
respectively.  Similarly, c switches with d, and e with f.  We can 
restore the square ice configuration by directing edges based on their 
color and the parities of their vertices; therefore, blue-degree $2$, 
green-degree $2$ graphs with the alternating boundary conditions are 
in bijection with square ice configurations with the vertical-out, 
horizontal-in boundary conditions, and thus with ASMs.  
Figure~\ref{same} shows an ASM represented as a matrix, as a 
square ice configuration, and as a colored graph.

In the subgraph of blue edges, all interior vertices have degree $2$,
and the endpoints have lower degree; thus, the connected components of
this graph are cycles, paths, and isolated points. The isolated points
are the green endpoints and may be ignored at the moment. The blue
endpoints are the only vertices of degree $1$, so all blue paths begin
and end at two blue endpoints. We call the two endpoints of a
monochromatic path {\em paired}\/, and refer to the partition of the
endpoints into these pairs as the {\em pairing}\/ of that ASM.

\section{The Theorem}
\label{theorem}
\begin{figure}[htb]
\begin{center}
\epsfig{file=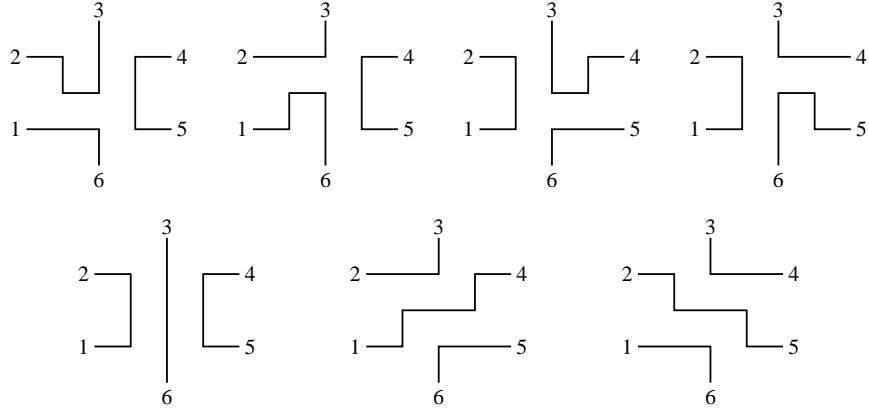,height=2.2in}
\end{center}
\caption{The blue edges in the seven ASMs of order $3$.}
\label{order3}
\end{figure}

Carl Bosley and Lukasz Fidkowski \cite{BF} conjectured
a general principal, which we illustrate for order $3$.
In the seven ASMs of order $3$, shown in 
Figure~\ref{order3}, each vertex is paired in three cases with its 
neighbor to the left, in three cases with its neighbor to the right, 
and once with the opposite vertex.  Let us number the blue endpoints 
clockwise, starting with $(0,1)$.  Then the observation is true for 
both vertex $1$, which is on a corner, and vertex $3$, which is in the 
middle of a side.  In general, their conjecture was that 
from the perspective of the pairing data, the $2n$ 
blue endpoints are arranged not around a square, but at the vertices 
of a $2n$-gon: if they are rearranged by an element of $D_{2n}$, the 
number of ASMs pairing blue endpoints $i$ and $j$ remains the same.  
Henry Cohn and James Propp \cite{CP} refined this 
conjecture to the form that we prove. This theorem, the central result 
of this paper, is stated below.

\begin{Theorem}
Let $\A_n(\pi_B,\pi_G,\ell)$ be the set of ASMs of order $n$ in which the 
blue subgraph induces pairing $\pi_B$, the green subgraph induces 
pairing $\pi_G$, and the sum of the number of cycles in the two 
subgraphs is $\ell$.  If $\pi_B'$ is $\pi_B$ rotated {\em clockwise}\/, 
and $\pi_G'$ is $\pi_G$ rotated {\em counterclockwise}, then  the sets 
$\A_n(\pi_B,\pi_G,\ell)$ and $\A_n(\pi_B',\pi_G',\ell)$ are in bijection.
\end{Theorem}

The bijection of the theorem is given by gyration, which is defined in 
the next section.  Since the blue and green endpoints rotate in 
opposite directions, we number the green endpoints clockwise by 
reflecting the blue labels over the line $y=x$.  This labeling 
compensates for the opposite behavior of the two colors, so that the 
effect of $G$ is to induce the same permutation of the labels: an 
increase by $1$.  In Section~\ref{gyration2}, 
we factor $G$ into two involutions, which have the same 
effect on the pairings as the generating reflections of $D_{2n}$.  
Thus if $\sigma\in D_{2n}$ is a permutation of the numbers from $1$ to 
$2n$ that induces $\pi_B'$ by rearranging the labels from the pairing 
$\pi_B$ and similarly induces $\pi_G'$ from $\pi_G'$, then 
$|\A_n(\pi_B,\pi_G,\ell)|=|\A_n(\pi_B',\pi_G',\ell)|$.

\section{Definition of Gyration}
\label{gyration}
Recall that we represent an ASM as a coloring of the graph $L_n$, which
lies in the square grid. The square grid is a planar graph and all of
its faces are squares.  Let us call the square with lower left corner
$(i,j)$ {\em even}\/ or {\em odd}\/ according to the parity of $i+j$.
Every edge in the grid is in one even and one odd square.  To each
square $S$ intersecting $L_n$, assign a function $G_S\:\A_n\to\A_n$.  In
the height function representation of an ASM The function is given by
fixing the color of all edges not in $S$ and determining the colors of
the four edges in that square based solely on their original colors.
If $S$ is on the boundary so that one or two of its edges are not in
$L_n$, let $G_S$ be the identity.  These boundary functions are not
strictly necessary, but because they ensure that the even (or odd)
squares that have functions cover $L_n$, they are a notational
convenience in the next section.  Otherwise, there are $2^4$, or $16$,
possibilities for the colors of the four edges.  In $14$ cases, two
edges in the square incident to the same vertex are the same color.
ASMs that have $S$ colored in one of these ways are fixed points of
$G_S$. In the remaining cases, the four edges of $S$ are alternately
blue and green; the two horizontal edges are one color, and the two
vertical the other.  For these two cases, the function reverses the
color of all four edges.  Since the only changes that $G_S$ makes is
to switch these two cases, it is an involution.

In the square ice view of these objects, reversing the colors 
corresponds to reversing the direction of oriented edges.  The 
involution reverses all four edges if and only if they are directed 
all clockwise or all counterclockwise about the square.  In the ASM 
view, the involutions still have a local effect, switching between 
$0$s and $\pm1$s, but the particular change and the decision to change
depend on entries 
other than the four corresponding to the vertices of the square.  
The simplest description of $G_S$ is in terms of the height function 
representation, which is a matrix with entries
offset from those in the ASM. 
Thus each entry lies on a square $S$ intersecting $L_n$. 
The matrix produced by $G_S$ agrees with the original matrix in every
entry except possibly the one on $S$. There are at most two possible
values that the entry on $S$ can take. If there is only one $G_S$
leaves the matrix unchanged. If there are two, $G_S$ changes the entry
to be the other possibility. 

Each $G_S$ may be considered to be ``local'' because it depends on and
affects only a small set of edges.  If $S$ and $S'$ are distinct even 
(or  odd) squares, then their edge sets are disjoint; moreover, 
$G_S$ and $G_{S'}$ commute. Thus we 
may, without worry about order, define $G_0$ as the composition of the 
involutions of even squares, and $G_1$ as the composition 
of involutions of odd squares. Since they are  compositions 
of commuting involutions, $G_0$ and $G_1$ are again involutions.
Finally, we define $G$ by $G_0\circ G_1$.

Let's summarize the definition: to perform a gyration on a graph,
visit each unit square in the graph $L_n$, first the odd ones and then 
the even ones.  In each visit, leave alone a square colored almost any 
way, but reverse the colors of all four edges if there are two 
parallel blue edges and two parallel green.

\section{Alternative Decomposition}
\label{gyration2}
It will be convenient to form another decomposition 
$G=H_0\circ H_1$, where the $H_k$ 
are no longer functions $\A_n\to\A_n$, but affect the paths much like $G$.
The range of $H_1$ and the domain of $H_0$ 
are no longer $\A_n$,
but instead the similar set $\B_n$ of graphs with the 
same underlying graph and the
same restriction on colors at vertices of degree $4$, but
color-reversed boundary conditions.
The blue endpoints alternate with the 
green ones, so gyration moves a blue endpoint two steps along, to the 
next blue one.
By stopping between $H_1$ and $H_0$, we will find an 
intermediate place, where the blue path has moved halfway to its 
destination and ends at what we have labeled a green endpoint.

Because of the local nature of its definitions, $G_S$ may be given a 
larger domain and range, such as $\A_n\cup\B_n$ (or even all 
colorings).  Then $G_0$ and $G_1$, defined by compositions, may also 
be extended to that domain.  Let $R$ reverse the color of each edge.  
Since $G_S$ is symmetric in blue and green, it commutes with the $R$.  
Defined by composition, $G_0$ and $G_1$ also commute with $R$.  Then 
define $H_k$, another involution, as the composition of $G_k$ with 
$R$. As $R^2$ is the identity, $H_0\circ H_1=G_0\circ G_1=G$.

\begin{figure}[htb]
\begin{center}
\epsfig{file=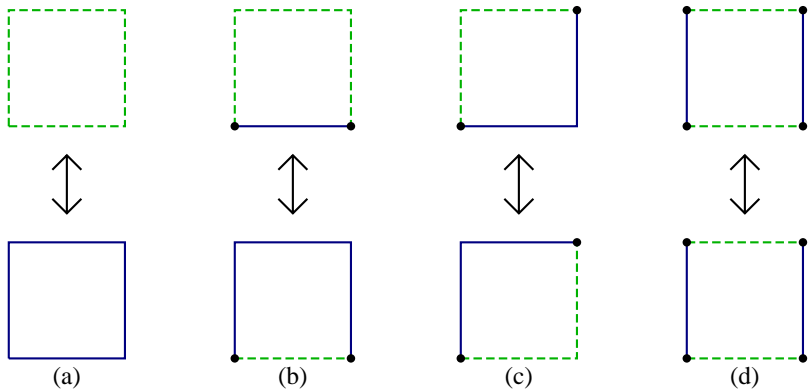,height=2.2in}
\end{center}
\caption{Examples of the effect of $H_S$ on the edges of $S$. 
Large vertices are {\em fixed\/} (see Section~\ref{lemma}) 
with respect to the parity of the squares shown.}
\label{local}
\end{figure}

Much as $G_0$ (respectively $G_1$) is the composition of local 
involutions associated with even (odd) squares, $R$ may be broken into 
the local reversals of even (odd) squares.  Define $H_S$ 
as $G_S$ composed with the reversal of the colors of the edges of $S$, 
and define $H_0$ (respectively $H_1$) as the composition of the (commuting) 
$H_S$, for $S$ even (odd).  Then $H_S$ preserves $2$ of the colorings of 
the edges of $S$ and reverses the colors in the other $14$.  
Figure~\ref{local} shows the effect of $H_S$ on $S$ for a complete 
(up to rotations) set of colorings of $S$.  

From the perspective of this decomposition, 
the summary of $G$ becomes the following: to perform gyration, 
visit each unit square, in the same order as the previous description, 
and reverse the colors of the edges, if any two of the same color are 
incident to the same vertex, or some edge of the square is missing 
from $L_n$; if, instead, the four edges alternate blue and green around 
the square, then leave them unchanged.

\section{The Bijection of Components}
\label{lemma}
Now that we have defined the $H_k$, we may begin to prove their 
properties.  Fix a graph $c\in\A_n$ and $k\in\{0,1\}$. Call an interior 
vertex {\em fixed}\/ if the unit squares of parity $k$ containing 
its two incident blue (or, equivalently, green) edges are distinct.  
The name is motivated by lemma below.
The fixed vertices are the points where a monochromatic path moves 
from one square of parity $k$ to another.
Fixed vertices provide the first suggestion that paths are well 
behaved under the application of $H_k$: a vertex is fixed in $c$ if 
and only if it is fixed in $H_k(c)$, as shown in Figure~\ref{local}.
This preservation also suggests the utility of fixed vertices; their 
drawback is that $H_{1-k}$ does not preserve them.
Instead of focusing on the endpoints, which $H_k$ changes, it is 
easier to work with the fixed vertices, which stay the same, as is 
shown in the following lemma, which is the principal step in the 
proof of the theorem.

\begin{Lemma} Two fixed vertices are in the same component of the 
blue (respectively green)
subgraph after the 
action of $H_k$ if and only if they were in the same component of the 
blue (green) before.
\end{Lemma}

\begin{proof}
First consider distinct fixed vertices connected by a path of nonfixed 
vertices.  Since fixed vertices are the points where paths move 
between squares of parity $k$, such a path must be contained in a 
single unit square of parity $k$; moreover, for a unit square to 
contain two fixed vertices, it must be contained in $L_n$ instead of
being on the boundary.
Figure~\ref{local} shows that, if 
there is a monochromatic path between two fixed vertices before the 
application of $H_k$, then there is one afterwards.  Specifically, in 
two cases (Figure~\ref{local}(a)), all four edges are the same color 
and there are no fixed vertices. In two cases 
(Figure~\ref{local}(d) and a rotation), 
the edges alternate colors and there are four fixed vertices. In 
these cases, $H_k$ does not change the colors, so we may reuse the 
old path. In the remaining twelve cases (Figures~\ref{local}(b),~(c) and 
rotations), there are two fixed vertices connected by one blue and 
one green path. The application of $H_k$ switches these colors, 
leaving one path of each color.

As the endpoints have degree $0$ or $1$, they cannot be in paths 
connecting interior vertices, such as fixed vertices.  Given a 
path between fixed vertices $v$ and $w$, we may 
divide it at each fixed vertex to obtain many paths of nonfixed 
vertices connecting fixed vertices.  From such paths we may 
construct paths of nonfixed vertices connecting the intermediate 
fixed vertices after the application of $H_k$. Concatenating these 
paths, we produce a path connecting $v$ and $w$.

Since $H_k$ is an involution, the lemma in one direction implies the 
converse.
\end{proof}

The lemma allows us to identify the paths before and after the application 
of $H_k$. It proves that paths behave well
and enables us to ask how they move.
The involution changes the path, except at the fixed
vertices, hence their name.
The theorem makes more sense from this perspective: the paths remain 
intact, but their endpoints circulate, changing the pairing.

Fixed vertices record where paths moved between squares of a 
particular parity. Since the 
application of $H_k$ does not change whether a vertex is fixed, 
nor, as the lemma tells us, which path passes through that vertex, a 
path must pass through the same sequence of squares of parity $k$ before 
and after the application of $H_k$. All that may change is how the 
path moves inside each of those squares and what it does before the first 
fixed vertex and after the last. All that remains of the proof of 
the theorem is to see 
what happens in the initial and final segments between the endpoints
and the fixed vertices.

\section{Proof of the Theorem}
\label{proof}
Take a graph $c\in\A_n(\pi_B,\pi_G,\ell)$.

As the blue subgraph has maximum degree $2$, its components are 
paths, cycles, and isolated vertices.  The lemma gives us a bijection 
between those blue components with fixed vertices before and after 
the application of $H_k$.  
A square of parity $k$, but not wholly 
contained in $L_n$, contains a blue endpoint, a green endpoint, a 
fixed vertex, and at most one nonfixed vertex.  The fixed 
vertex is connected, by appropriately colored paths, to both 
endpoints.  Thus every path connecting two blue endpoints contains at 
least one fixed vertex.  
Moreover, this fixed vertex is connected to the endpoints after 
$H_k$, so the bijection of components turns paths into paths and not 
into cycles.

Thus we have a bijection between the blue paths connecting endpoints
before the application of $H_k$ with the blue paths connecting endpoints 
afterwards. 
We may compose the bijection induced by applying $H_1$ to 
$c$ with that induced by applying $H_0$ to $H_1(c)$ to obtain a 
bijection of the sets of blue paths in $c$ and $G(c)$. Since paths 
induce pairings of the endpoints, and we can associate a path in $c$ 
with a path in $G(c)$, we can associate the corresponding pairings of 
the endpoints of the two paths. All we need is to determine how the 
end of a path moves.

Recall that if a unit square of parity $k$ contains an endpoint, it is on the
boundary of $L_n$ and it contains a fixed vertex and two endpoints,
one blue and one green.
When one of the $H_k$ affects
that square, it switches the endpoints of the blue and green paths,
leaving them both passing through the fixed vertex.  For example, the
unit square indexed by $(0,1)$ and affected by $H_1$ contains the blue
endpoint $(0,1)$ and the green endpoint $(0,2)$.  One of the vertices
$(1,1)$ and $(1,2)$ is fixed, depending on whether the edge between
them is green or blue, respectively.  The effect of $H_1$ is to switch
the endpoints, to make a blue path end at $(0,2)$ and a green path end
at $(0,1)$.  The endpoint of the blue path travels clockwise and the
endpoint of the green path counterclockwise.  Since there are two
endpoints in each unit square on the boundary and the endpoints
alternate colors, in each unit square the end of the blue path starts
counterclockwise of the end of the green path and ends clockwise of
it.  Similarly, the unit squares associated with $H_0$ contain one
endpoint of each color.  They are arranged in the other order, with
what is normally a blue endpoint clockwise of what is normally a green
endpoint, but $H_1$ has already reversed the colors before the
application of $H_0$.  Thus $H_0$ also sends the blue paths clockwise
and the green paths counterclockwise.

Over the course of the two steps of $G$, the end of a blue path moves 
two endpoints, or one blue endpoint, clockwise.  Thus the label on 
each endpoint of a blue path increases by $1$. The effect on the 
pairing is the same as if we had relabeled the endpoints in that way.
Since the endpoints move by switching blue and green, the green 
ones move counterclockwise, as was anticipated in the choice of green labels.

Finally, we must determine the number of closed loops of both paths in 
$G(c)$.  The lemma states that $H_k$ turns cycles that 
contain fixed vertices into new cycles of the same color that 
contain fixed vertices; so the number of cycles of each color 
containing a fixed vertex is preserved.  A cycle that does not 
contain a fixed 
vertex must be contained in a single unit square of parity $k$.  
Such a cycle reverses color; the total number of such cycles of both 
colors is preserved. Thus each $H_k$ preserves the total number of 
cycles; so must $G$, which is $H_0\circ H_1$. In addition, this 
preservation is bijective, although the division into cycles with 
fixed vertices and those without is not preserved because it depends 
on $k$.

The three statistics match the claims of the theorem, so we have the 
inclusion,
$$G(\A(\pi_B,\pi_G,\ell))\subset
\A(\pi_B',\pi_G',\ell).$$ 
The reverse inclusion can be obtained by a similar argument about 
$G^{-1}$ because $H_0$ and $H_1$ played essentially identical roles 
in the proof and $G^{-1}=H_1\circ H_0$.

To establish the dihedral version of the theorem, note that the 
function $d$ that reflects a graph across the line $y=x$ reverses the 
colors of the endpoints.  Since $H_1$ and $H_0$ also reverse the 
colors of the endpoints, composing $d$ with either of them gives a 
function $\A_n\to\A_n$.  Moreover, the reflection preserves the even 
and odd squares, so it commutes with $H_1$ and $H_0$.  Thus $H_kd$ is 
an involution, and $G=H_1d\circ H_0d$.  The effect of $H_1d$ is to 
send Endpoint $i$ to Endpoint $2n+1-i$, and $H_0d$ sends $i$ to 
$2n+2-i$ (modulo $2n$); in particular, $H_0d$ fixes $1$ and $n+1$.  These 
two reflections generate $D_{2n}$.  \hfill\qedsymbol

\section{Generalizations and Conjectures}
\label{generalizations}
If $n<5$, then $G^{2n}$ is the identity. For larger $n$, it is a 
bijection $\A_n\to\A_n$, which preserves the three statistics, but which 
is not the identity. Its order is unknown in general. 
Table~\ref{order} shows that the order tends to be divisible by rather 
small primes, but this phenomenon is not surprising as the order is the 
least common multiple of the sizes of the orbits. A similar function 
is $G^n$ composed with rotation by $\pi$. 

\begin{table}[h]
\begin{center}
\caption{The order of $G^{2n}$}
\label{order}
$$\begin{array}{cc}
\hline\noalign{\smallskip}
n & \textrm{order}\\
\hline\noalign{\medskip}
4 & 1\\
5 & 2\\
6 & 2\cdot3\cdot5\cdot7\\
7 & 2^2\cdot3^2\cdot5\cdot7\cdot11\cdot19\\
\hline\noalign{\medskip}
\end{array}$$
\end{center}
\end{table}

A torus graph with both dimensions even is bipartite, so the bijection 
between square ice and colored graphs remains.  Such a torus graph has 
bipartite dual, allowing the division of the squares into even and 
odd, which allows the two steps of gyration.  Gyration preserves the 
total number of cycles, but there are no paths to change.  Much as 
with $G^{2n}$ on the square, we ask for the order of gyration on the 
torus.  Similarly, gyration is well defined and nicely behaved on the
infinite plane, although counting is more difficult and statistics may
not be defined.

Bijections similar to gyration may be constructed on many models that
have a height function, though the representation in terms of colored
paths is lost and with it the theorem.  Mills, Robbins, and Rumsey
examined involutions analogous to $G_0$ and $G_1$ under
the names $\rho$ and $\gamma$ \cite{MRR83}.

Knowing that certain sets are in bijection, we are led to ask for 
their cardinalities.  Let $A_{n,k}$ be the number of ASMs 
of order $n$ such that the blue paths pair vertex $i$ with vertex 
$2k+1-i$ for $1\le i\le 2k$, so that there are nested paths $k$ deep.  
Let $B_{n,k}$ be those ASMs with that constraint and the additional 
requirement that for $i>k$, blue endpoint $2i$ is paired with blue 
endpoint $2i-1$.  
The author conjectures 
that $A_{n,k}=B_{n+1,k+1}$.  Equality holds for $1\le n\le 8$ and $0\le 
k\le 5$.  In the cases $k=0,1$, the conjectures are due to
James Propp and David Wilson \cite{PW}.

The case $k=0$ gives $A_{n,0}=B_{n+1,1}$, which is the most 
interesting case because both sides simplify.  The set counted by the 
right-hand side is the same as that counted by $B_{n+1,0}$ because 
they both require that vertex $1$ pair with vertex $2$. The new name 
$B_{n+1,0}$ is more appealing than $B_{n+1,1}$ because it unifies 
all of the pairings under a fairly simple rule: every blue path 
connects consecutively  numbered endpoints. The conjecture implies 
that this set is enumerated by $A_{n,0}$, which is the number of 
unrestricted ASMs, and is shown in $\cite{Z}$ and $\cite{K}$ to be 
$\prod_{i=0}^{n-1}\frac{(3i+1)!}{(n+i)!}$.

\end{document}